\documentclass[11pt]{article}
\usepackage{amscd,amssymb,stmaryrd}
\usepackage{amsmath}
\usepackage{graphicx}
\usepackage{subcaption}
\usepackage[utf8]{inputenc}
\usepackage[export]{adjustbox}
\usepackage{wrapfig}
\usepackage{amstext}
\usepackage{pstricks,pst-node,pst-plot,pst-coil}
\usepackage{amsthm}
\usepackage{amsmath}
\usepackage{color}
\usepackage{mathtools}
\usepackage{hyperref}
\usepackage[english]{babel}
\usepackage{booktabs}
\usepackage{mathrsfs}
\usepackage{comment}
\usepackage{float}
\usepackage[linesnumbered,ruled,vlined, ruled,vlined]{algorithm2e}

\theoremstyle{definition}

\theoremstyle{remark}
\newtheorem*{remark}{Remark}

\title{Neural homogenization and the physics-informed neural network for the multiscale problems }
\author{Wing Tat Leung\footnote{Department of Mathematics, University of California, Irvine, USA}, ~ Guang Lin\footnote{Department of Mathematics and Mechanical Engineering, Purdue University}, ~ Zecheng Zhang \footnote{Department of Mathematics, Purdue University} }
\begin{document}
\maketitle

\begin{abstract}
  Physics-informed neural network (PINN) is a data-driven approach to solve equations. 
    It is successful in many applications;
    however, the accuracy of the PINN is not satisfactory when it is used to solve multiscale equations.
    Homogenization is a way of approximating a multiscale equation by a homogenized equation without multiscale property;
    it includes solving cell problems and the homogenized equation.
    The cell problems are periodic; and we propose an oversampling strategy which greatly improves the PINN accuracy on periodic problems.
    The homogenized equation has constant or slow dependency coefficient and can also be solved by PINN accurately.
    We hence proposed a 3-step method, neural homogenization based PINN (NH-PINN), to improve the PINN accuracy for solving multiscale problems with the help of the homogenization. 
    We apply our method to solve three equations which represent three different homogenization. 
    The results show that the proposed method greatly improves the PINN accuracy in particular when the scaling is small.
    Besides, we also find that the PINN aided homogenization may achieve better accuracy than the numerical methods driven homogenization; PINN hence is a potential alternative to implementing the homogenization.
\end{abstract}

\section{Introduction}
In the recent decade, the problems having features at multiple scales are of rising importance. For example, in biomedical applications, to obtain a more accurate model, people need to incorporate the cell-level information to the tissue-level model. In material application, people are interested in the macroscopic properties of some complex composite materials which are composited in a microscopic scale. In petroleum application \cite{popov2009multiphysics, ma2008efficient, efendiev2000modeling, gildin2013nonlinear}, the fluid flow in the reservoir is heavily depending on the heterogeneous of the porous media.
\textcolor{black}{
It is not easy to solve the equations by the classical methods.
To capture the multiscale properties, one needs to use a fine mesh \cite{chung2016adaptive, chung2018constraint, efendiev2000convergence, chung2021computational, chetverushkin2021computational}; or one can use the multiscale finite element methods which can save the computation cost but it is still not easy for some problems.
Hence researchers also studied the data driven approaches to solve the multiscale problems.
}

Physics-informed neural network (PINN) is a neural network approach of solving the partial differential equations (PDE) \cite{raissi2019physics, liu2021deep, lin2021multi}. 
The idea of the PINN is to approximate the solution of the PDE by a network.
The PINN has been widely used in solving both forward and inverse problems \cite{yang2021b, lin2021multi}.
Compared to the classical numerical methods, PINN is a mesh free method and hence can interpret the solution without a mesh.
Moreover, PINN has no CFL constrain and is easy to be used to solve time-dependent problems.
For the equation with a convection term, PINN is also powerful and one does not need to consider the directions of the flow, 
which reduce the difficulty. Besides, PINN can be used to solve the equations with the mixed direction derivatives such as paraxial approximation \cite{chung2021computational}. The traditional numerical methods are usually hard to implement and involve intense computation.
Compared to other neural network methods of solving the PDE \cite{zhang2020learning, chung2021multi}, PINN is a sample free method which does not rely on label,
we hence can use PINN on solving problems not limited to UQ.

Although researchers have applied PINN to solve many mathematical problems, there are only a few works about multiscale problems.
In \cite{wang2021eigenvector}, the authors point out that the classical PINN formulation is unable to capture the multiscale property of the solution; 
they later proposed a method which applies the Fourier transformation to prepossess the input to solve the problem.
Their interest in the multiscale is in the solution; this is different from many existing real life multiscale problems we have discussed before \cite{hou1997multiscale, efendiev2013generalized, zhang2020learning, chung2018cluster, chung2020multi}. \textcolor{black}{In this work, we mainly focus on the multiscale problems whose multiscale property comes from the equation. 
One typical example is the 2D elliptic equation which models the porous media flow. 
The permeability of the equation has a multiscale nature and brings multiscale to the solution; for example,
the permeability $\kappa(x) = \sin(x)+\sin(10x)$.
}

We hence seek help from the PINN method; however, we find that the classical PINN cannot give us an accurate prediction.
The relative error is very large and the training is not robust to the hyper-parameters and randomness. 
More precisely, the relative error varies a lot for different runs even if we keep all hyper-parameters the same.
To alleviate the problems of the classical PINN and take the advantages of the PINN, we propose a homogenization based approach.

Homogenization is one of the classic ways for solving an equation with highly heterogeneous medium \cite{bensoussan2011asymptotic, engquist2008asymptotic, pankov2013g, tartar1979compensated}. Basically, throughout the homogenization process, we will obtain a homogenized equation of the original equation such that the solution of this homogenized equation can capture the macroscopic behavior of the solution of the original equation. The parameters of the homogenized equation are normally homogeneous or smooth such that the equation can be solved numerical in a coarse mesh. For example, we consider a diffusion equation with a periodic media $\kappa$:
\begin{equation}
-\nabla \cdot (\kappa(\cfrac{x}{\epsilon})\nabla u_{\epsilon} ) = f.
\end{equation}
By homogenization theory, we can obtain a homogenized equation: 
\begin{equation}
-\nabla \cdot (\kappa^*\nabla u_{0} ) = f.
\end{equation}
where $\kappa^*$ is the homogenized parameter which is a constant tensor, and we have $u_0$ is converging to $u_{\epsilon}$ in $L_2$-norm as $\epsilon\rightarrow 0$. Thus, one can solve the homogenized equation in a coarse grid and use the homogenized solution $u_0$ to approximate the solution of the original equation $u_{\epsilon}$.

We now propose to solve the multiscale problems by PINN with the help of homogenization \cite{efendiev2009multiscale}.
We find that the homogenization can be implemented by PINN very naturally.
The classical homogenization consists of three steps. The first step is to solve the cell problems on a unit cube. 
The cell problems are equipped with a periodic boundary condition and the permeability has no fast dependency.
It should be noted that we find the periodic PDE is not easily solved by PINN when the dimension is high; 
however, we propose an oversampling trick which greatly improves the performance of solving the periodic problems.
In all, the cell problems can be solved by PINN method accurately and efficiently.

The second step is to evaluate the homogenized permeability which has no multiscale property and then obtain the homogenized equations.
This step usually involves finding the derivatives of
the cell problem solutions. 
If the solutions are approximated by a network, the derivatives can be easily derived by the auto differentiation of the software.
Once the cell problems are solved by the network accurately, we can anticipate getting an accurate homogenized equation.

Lastly, we can solve the homogenized equation. 
The permeability associated with this equation is usually constant or depends on the slow variable only,
the PINN has shown its power in solving such kind of equations.
We conclude that the PINN can be used to implement the entire homogenization process; 
one is able to solve the multiscale PDEs by PINN with the help of the homogenization.

We summarize the contributions of this work as follows:
\begin{enumerate}
     \item We observe that the accuracy of the PINN degenerates when solving the multiscale PDE; \textcolor{black}{we also explain the reason of the failure.} 
    Please note that in this work, we focus on the PDE whose multiscale property originates from the equation coefficients,
    in particular, we are interested in the multiscale permeability with multiple frequencies.
    \item We propose a 3-step approach, neural homogenization based PINN (NH-PINN). Instead of solving the multiscale equation by PINN directly, we first apply PINN to solve the cell problems which are used to derive the homogenized equation; the homogenized equation can then be easily solved by PINN.
    \item We propose an oversampling strategy to solve the periodic PDE by PINN; this method greatly improves the accuracy of the PINN when solving the high dimensional periodic problems.
    \item We conduct four numerical experiments which represent three different types of homogenization.
    The predictions are very accurate and are much better than the classical PINN in particular when the scaling is small.
    \item In addition to the improved accuracy of the PINN, we also observe that NH-PINN can improve the homogenization accuracy.
    That is, if we apply PINN to implement homogenization, the solution may be more accurate than the traditional numerical methods (eg, finite element methods) driven homogenization. \textcolor{black}{We hence suggest that PINN may be a potential alternative of implementing the homogenization.}
\end{enumerate}

The rests of the paper are organized as follow. 
In Section \ref{pinn_sec}, we review the basics of the PINN; the performance of applying classical PINN to solve multiscale problems is also presented in this Section. Next, we will review the homogenization method in Section \ref{homog_sec}. We will give the details of one homogenization of the elliptic equations. Our numerical examples are not limited to the elliptic operator, we hence also present the homogenization of the reaction diffusion equation in the Appendix (\ref{app_sec}). We conduct four experiments and they are shown in the Section \ref{numerical_sec}.

\section{Physics-informed neural network (PINN)}
\label{pinn_sec}
In this section, we briefly review the physics-informed neural network (PINN) for solving the PDE; for the inversion problems, please refer to \cite{yang2021b, lin2021multi} for details. 
PINN is a data driven method to solve the PDE.
The solution of the PDE is approximated by the network and the target is to minimize the error of this approximation.
The PINN has been applied to solve various of problems; however, it is rarely used in solving multiscale problems.
Suppose we are solving the following system in the domain $\Omega$,
 \begin{align*}
     \mathcal{L}(u) &= f \text{ in } \Omega \\
     \mathcal{B}(u) &= b \text{ on } \partial \Omega
 \end{align*}
 where $\mathcal{L}$ is a differential operator and $\mathcal{B}$ is the boundary condition operator. $f$ is the given source term, $b$ is the given boundary condition. 
 The idea of the PINN is to build a map from $\Omega$ to the solution by a network $ \mathcal{F}_{\beta}(\cdot)$, where $\beta$ is the parameters associated with the network. More precisely, we will solve the following minimization problem:
 \begin{equation}
     \min_{\beta} \frac{w_1}{N_f} \sum_{i=1}^{N_f} |\mathcal{L}(\mathcal{F}_{\beta}(p_i)) - f(p_i)|^2 + \frac{w_2}{N_b} \sum_{i=1}^{N_b} |\mathcal{B}(\mathcal{F}_{\beta}(q_i)) - b(q_i)|^2,
     \label{pinn_pinn}
 \end{equation}
 where $w_1+w_2 = 1$ are the positive weights; 
 $\{p_i\}\subset\Omega$, $\{q_i\}\subset\partial \Omega$ and
 $N_f, N_b$ are the number of points used in discretizing the domain and boundary respectively. 
It should be noted that if the PDE has other constrains such as the initial condition $\mathcal{I}$, we just need to sample some points $z_i$,
substitute the points in the network $\mathcal{F}_{\beta}(\cdot)$, then approximate the constrains by $\mathcal{I}(\mathcal{F}_{\beta}(\cdot))$
and include this loss $\mathcal{I}(\mathcal{F}_{\beta}(\cdot))- I(\cdot)$ in Equation (\ref{pinn_pinn}), where $I(\cdot)$ is the given innitial conditions.
PINN is easy to be implemented. 

Since $\mathcal{F}_{\beta}$ approximate the solution and is smooth (neural network), the differential operator $\mathcal{L}(\mathcal{F}_{\beta}(\cdot))$ can be evaluated with the auto differentiation package of the modern deep learning software easily. The minimization problem can then be solved by the variants of the gradient descent algorithm.
 
 \subsection{Failure of the PINN}
\begin{figure}[H]
\centering
\includegraphics[scale = 0.3]{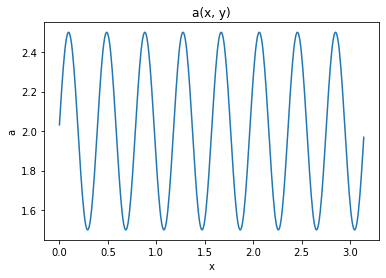}
\includegraphics[scale = 0.3]{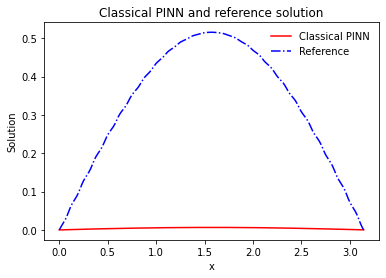}
\includegraphics[scale = 0.3]{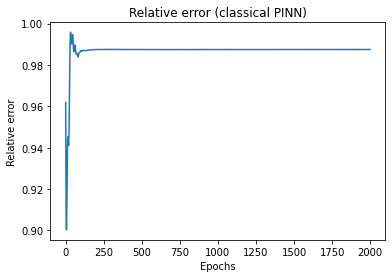}
\caption{1D elliptic problem without slow dependency (classical PINN). Left: demonstration of the permeability $\kappa(x) = 0.5\sin(2\pi x/\epsilon)+2$, where $\epsilon = 1/8$. Middle: learnt solution by the classical PINN vs the reference solution. Right: relative error as a function of the training epochs. The average error of the last 500 epochs is $0.987471$.
In this example, learning rate is 0.0001 and $w_1 = w_2 = 0.5$ and we train the network with Adam gradient descent for 2000 epochs.
We test with different sets of the hyper-parameters such as the learning rate and the weights of the loss, however, we cannot obtain a satisfied result.
} 
\label{pinn_1d_fast}
\end{figure}

 PINN has been used to solve various problems, however the PINN performance of solving multiscale problems has been compromised.
 For example, let us consider the following elliptic equation:
 \begin{align*}
     & -\nabla \cdot (\kappa \nabla u(x)) = f, x\in \Omega,\\
     & u = 0, x\in\partial \Omega,
 \end{align*}
where $\Omega = [0, \pi]$ and $f = \sin(x)$.
This is a standard elliptic problem with a homogeneous Dirichlet boundary condition. If $\kappa(x)$ is a constant or has no multiscale property, the PINN performs very well; however if
the permeability $\kappa(x)$ presents some multiscale property, for example, $\kappa(x) = 0.5\sin(2\pi x/\epsilon)+2$, where $\epsilon = \frac{1}{8}$, we cannot obtain a satisfied solution (with a low relative error) using the classical PINN.
We solve the above problem and demonstrate the results in 
Figure (\ref{pinn_1d_fast}). 
We will give more examples about the direct application of PINN on solving multiscale problems later in Section (\ref{numerical_sec})
and provide our homogenization strategy; that is, we propose to solve the multiscale problems in three steps with the help of the homogenization.
Our numerical experiments show that the accuracy of the PINN is greatly improved.

\subsection{Intuitive explanation of the failure}
\label{pinn_reason_sec}
Now the question is why PINN fails in solving the multiscale equation.
We give an explanation motivated by the finite element methods (FEM).
For the standard FEM,  it is well known that a coarse mesh FEM solution cannot resolve the multiscale solution. 
For example, in 1-d case, the multiscale solution $u_{\epsilon}$ satisfying 
\[
-\partial_x \Big{(}a(\frac{x}{\epsilon}) \partial_x u_{\epsilon}\Big{)} = f,
\]
converging to the homogeneous solution $u_0$ satisfying
\[
-\partial_x \Big{(}a^*(x) \partial_x u_{0}\Big{)} = f,
\]
 in $L_2$-norm when $\epsilon$ goes to $0$ where $a^*$ is the harmonic average of $a$ \cite{efendiev2009multiscale}.
However, if we use linear finite element method in a coarse grid, we have $\partial_x \phi_i$ is a piecewise constant for any linear basis function $\phi_i$. Thus, we have the finite element solution $u_{FEM}$ satisfying 
\[
\int_{\Omega} a(\frac{x}{\epsilon}) \partial_x u_{FEM}\partial_x \phi_i \approx \int_{\Omega} \overline{a}(x) \partial_x u_{FEM}\partial_x \phi_i, \forall \phi_i 
\]
where $\overline{a}$ is the average of $a$. We can check that the harmonic average and normal average can have a large difference in some cases. For example, if $a=15 \sin^2(x)+1$, we have $a^*=4$  but $ \overline{a}=8.5$. Thus, we have $u_{FEM}\approx \frac{8}{17}u_{0}\approx u_{\epsilon}$ and the relative error is $\cfrac{\|u_{FEM}-u_{\epsilon}\|_{L_2}}{\|u_{\epsilon}\|_{L_2}}\approx \cfrac{9}{17}>50\%$. 
We remark that using higher order coarse mesh finite element method in a higher dimensional case will have a similar behavior.

It is shown in paper \cite{rahaman2019spectral} that the network tends to learn the low frequency solution (globally varies without local fluctuation) first. The low frequency solution is similar to the solution in a coarse grid. If we use the multiscale loss function to train the solution, it will tend to obtain a FEM solution-like solution since when the solution is smooth, the multiscale residual is similar to the the residual with average parameters. PINN hence fails in solving multiscale problems.

One way to verify our hypothesis is: we can initialize the PINN network such that it can capture the low frequency component of the solution;
if the network cannot capture the low frequency as the training goes, 
this implies that the multiscale operator misleads the network and then the training fails.
Our proposed method (NH-PINN) can train a network which learns the coarse scale solution by the homogenization; we hence can use NH-PINN network as the initialization to the classical PINN network; and then train the network further with the classical PINN loss and settings. Our experiments show that the coarse scale solution initialization cannot improve the classical PINN. 
This indicates that the set of parameters which results in an approximation of the coarse scale solution may not be an optimum of the PINN loss.
The optimization algorithm performs well;
however the loss function itself cannot resolve the multiscale features.
This verifies our hypothesis of the PINN failure. 
The details of the experiments can be found in Section \ref{numerical_sec}.

\section{Homogenization}
\label{homog_sec}
As we mentioned in the introduction, to obtain the homogenizated parameter, we need to use the solution of the cell problems. For the completeness of the paper, in this section, we briefly revisit the derivation of the homogenizated equation of a diffusion problem and introduce the cell problem in this case. We also included a discussion about the homogenization of a diffusion-reaction equation in the appendix.

We here give an example of the homogenization of multi-dimension elliptic problems.
Consider we are solving:
\begin{align}
    -\frac{\partial }{\partial x_i}\bigg( 
    a_{ij}(x/\epsilon) \frac{\partial}{\partial x_j} u_{\epsilon}(x)
    \bigg) = f(x), x\in\Omega
\end{align}
with $u_{\epsilon}(x) = 0$ on $\partial\Omega$; here we use the Einstein notation.
We seek for $u_{\epsilon}(x)$ in the asymptotic expansion:
\begin{align}
    u_{\epsilon}(x) = u_0(x, x/\epsilon)+\epsilon u_1(x, x/\epsilon)+
    \epsilon^2 u_2 (x, x/\epsilon)...
    \label{homg_asym}
\end{align}
where $u_j(x, y)$ are periodic in $y = x/\epsilon$.
Denote 
\begin{align*}
    A^\epsilon = -\frac{\partial }{\partial x_i}\bigg( 
    a_{ij}(x/\epsilon) \frac{\partial}{\partial x_j}\bigg)
\end{align*}
It is not hard to check that:
\begin{align*}
    A^{\epsilon} = \epsilon^{-2} A_1+\epsilon ^{-1} A_2+ 
    \epsilon^0 A_3,
\end{align*}
where ,
\begin{align*}
    & A_1 = -\frac{\partial }{\partial y_i}\bigg( 
    a_{ij}(y) \frac{\partial}{\partial y_j}\bigg),\\
    & A_2 = -\frac{\partial }{\partial x_i}\bigg( 
    a_{ij}(y) \frac{\partial}{\partial y_j}\bigg)
    -\frac{\partial }{\partial y_i}\bigg( 
    a_{ij}(y) \frac{\partial}{\partial x_j}\bigg)\\
    & A_3= -\frac{\partial }{\partial x_i}\bigg( 
    a_{ij}(y) \frac{\partial}{\partial x_j}\bigg).
\end{align*}
We hence have $A_1 u_{\epsilon} +A_2 u_{\epsilon} + A_3 u_{\epsilon} = f$. Equating the terms with the same power, it follows that,
\begin{align}
    & A_1 u_0 = 0, \label{homog_u0} \\
    & A_1 u_1 + A_2 u_0 = 0, \label{homog_u1}\\
    & A_1u_2+A_2u_1+A_3u_0 = f. \label{homog_u2}
\end{align}
Substitute the $A_1$ and the equation (\ref{homog_u0}) becomes:
\begin{align*}
    -\frac{\partial }{\partial y_i}\bigg( 
    a_{ij}(y) \frac{\partial}{\partial y_j}\bigg) u_0 = 0
\end{align*}
The theory of the second-order ODE implies that $u_0$ is independent of $y$ and this will further simplify the equation (\ref{homog_u1}); it follows that,
\begin{align*}
    -\frac{\partial }{\partial y_i}\bigg( 
    a_{ij}(y) \frac{\partial}{\partial y_j}\bigg) u_1
    = \bigg( \frac{\partial }{\partial y_i}
    a_{ij}(y) \bigg) \frac{\partial}{\partial x_j} u_0.
\end{align*}
$u_1(x, y)$ can be solved by introducing $\chi_j(y)$ which is the solution of the problem:
\begin{align}
    &-\frac{\partial }{\partial y_i}\bigg( 
    a_{ij}(y) \frac{\partial}{\partial y_j}\bigg) \chi_j = 
    \frac{\partial}{\partial y_i}a_{ij}(y),\\
    & \chi_j \text{ is periodic in $y$ with mean $0$}.
    \label{homog_cell_prob}
\end{align}
The above problems (\ref{homog_cell_prob}) are called cell problems and they are to be solved in one period of $y$, or, in the unit cell $Y = [0, 1]^d$ where $d$ is the dimension of the problem.
$u_1$ can be then written as:
\begin{align*}
    u_1(x, y) = \chi_j \frac{\partial u_0}{\partial x_j}(x).
\end{align*}
Finally, we have,
\begin{align*}
    -\frac{\partial }{\partial y_i}\bigg( 
    a_{ij}(y) \frac{\partial}{\partial y_j}\bigg) u_2 =
    A_2u_1+A_3u_0-f.
\end{align*}
The ODE has a solution only if the right hand side has zero mean in one period of $y$, i.e., 
\begin{align*}
    \int_Y (A_2u_1+A_3u_0-f )dy = 0.
\end{align*}
This solvability condition then gives,
\begin{align}
    -\frac{\partial }{\partial x_i}\bigg( 
    a^*_{ij} \frac{\partial}{\partial x_j}\bigg) u_0 = f,
\end{align}
where $a^*_{ij} = \int_Y (a_{ij}+a_{ik}\frac{\chi_j}{\partial y_k}) dy$ is called the homogenized coefficient and $u_0$ is the homogenized solution.
\begin{remark}
\textcolor{black}{
Note that $a_{ij}(x/\epsilon)$ has no slow dependency $x$. If $a_{ij}(\cdot)$
also depends on slow variable $x$, the cell problems solutions also have the $x$ dependency, i.e., 
\begin{align*}
    -\frac{\partial }{\partial y_i}\bigg( 
    a_{ij}(x, y) \frac{\partial}{\partial y_j}\bigg) \chi_j(x, y) = 
    \frac{\partial}{\partial y_i}a_{ij}(x, y),
\end{align*}
$u_1(x, y)$ then becomes:
\begin{align*}
   u_1(x, y) =  \chi_j(x,y) \frac{\partial u_0}{\partial x_j}(x)
\end{align*}
The homogenized coefficient $a^*(x)$ then has slow dependency since we only integrate in $y$ of $\chi_j(x, y)$. This case will be illustrated by in Section (\ref{numerical_sec_1dslow}).}
\end{remark}

We now briefly discuss the convergence property of the homogenization. The goal is to estimate the remainder $R$ which is defined as:
\begin{align*}
    u_{\epsilon} = u_0 +\epsilon u_1+R.
\end{align*}
\textcolor{black}{$R$ then satisfies that $|R|<C\epsilon$ where $C$ is a constant. We also have the energy estimate $\int a |\frac{d}{dx} R|^2\leq C\epsilon$.}

\section{Neural homogenization based PINN (NH-PINN)}
\label{pp_sec}
In this section, we introduce our proposed method: neural homogenization based PINN (NH-PINN).
We have seen in Section \ref{pinn_sec} (and will see more examples in Section \ref{numerical_sec}) that PINN is unable to solve the multiscale problems with a satisfied accuracy; 
however from the Section \ref{homog_sec}, \textcolor{black}{the homogenized equation loses its multiscale property} since the homogenized coefficient is either constant or depends only on the slow variable;
besides, the previous studies of the PINN have shown that the PINN can deal with equations with constant coefficients; this motivates us to use homogenization.

The idea is to decompose the original hard problem into easier and solvable problems. Asymptotic expansion homogenization consists of three steps and
we have shown that the homogenized solution converges to the real solution in Section \ref{homog_sec};
if each of the steps can be solved by PINN easily, we believe that the accuracy of the final data driven (by PINN) solution can approximate the real solution closely.
Our method consists of three steps:
\begin{enumerate}
    \item Solve the cell problems using PINN.
    \item Evaluate the homogenized coefficients.
    \item Solve the homogenized equation using PINN.
\end{enumerate}
The first step is to solve the cell problems. 
Since there is no fast dependency, the cell problems can be solved easily by PINN.
It should be noted that the cell problems are equipped with the periodic boundary condition; 
we find that the performance is satisfactory when solving the 1D periodic problem; 
however, for the high dimensional cases, the PINN cannot give a good result. 
To deal with this issue, we propose an oversampling trick which facilities the PINN with a few more sampling points. 
Our numerical experiments show that this trick greatly improves the PINN performance for handling $2D$ periodic problem.
We introduce this trick in detail in Section (\ref{numerical_oversampling}).

The second step is to evaluate the coefficients. 
This step usually involves calculating the derivatives of the cell problem solutions.
The cell problems are solved by PINN, i.e., the solution is approximated by a network which is smooth, hence with the help of the auto differentiation, these derivatives can be calculated very easily.
This is one of the benefits of our method.

The last step is to solve the homogenized equation with PINN. Since PINN has shown its power in solving the equations with constant or slow varying coefficients, we do not expect any difficulty in this step. 

\textcolor{black}{There are three error sources in our problem. Firstly, the error comes from the homogenization, however, this error is small, in particular when $\epsilon$ is small and the scales are very different. The other two errors come from the PINN solver accuracy when solving the cell problems and the homogenized equation;
we will show in the numerical section that all three errors can be controlled but the error of solving the homogenized equation dominates the total error of our proposed method. }

\subsection{Oversampling}
\label{numerical_oversampling}
As we have discussed before, each cell problem is equipped with a periodic boundary condition.
We find that the generalization error is large when we solve the $2D$ equation with a periodic boundary condition. 
Hence, we are going to introduce an oversampling trick which is used to aid PINN to solve an equation with periodic boundary conditions.
Our experiments show that this trick greatly improves the training for the $2D$ problems.
For simplicity, will illustrate this idea in $1D$ case and the extension to the high dimensional cases are similar.

Suppose $\Omega = [x_0, x_t]$ and recall (\ref{pinn_pinn}), we then have $x_0 = q_1$ and $x_t = q_2$. Since we are solving the periodic problem,
the minimization problem then becomes:
\begin{align}
    \min_{\beta} \frac{w_1}{N_f} \sum_{i=1}^{N_f} |\mathcal{L}(\mathcal{F}_{\beta}(p_i)) - f(p_i)|^2 + \frac{w_2}{2} \sum_{i=1}^{2} |(\mathcal{F}_{\beta}(q_1) - \mathcal{F}_{\beta}(q_2)|^2,
     \label{numerical_pinn}
\end{align}
where the second term above indicates the periodic boundary condition. To enforce the boundary learning, we oversample some points on both sides of the domain. 
Denote $\Delta x = 1/(N_f-1)$, we then 
include additional 2 sets of points $\{q_i^{left}\}_{i = 1}^{N_o}$ and $\{q_i^{right}\}_{i = 1}^{N_o}$
, where $N_o$ is the number of oversampling layers and $N_o<N_f-2$.
In this work, we consider the uniform mesh, we hence have,
$q_{i+1}^{left}-q_{i}^{left} = \Delta x$ and $q_{i+1}^{right}-q_{i}^{right} = \Delta x$. We also require:
\begin{align*}
    &q^{left}_{N_o}+x_t-x_0 = x_t - \Delta x, \\
    &q^{right}_1-(x_t-x_0) = x_0+\Delta x.
\end{align*}
Due to the above assumptions, we have $q_i^{left}+x_t-x_0\in\{p_j\}$ and $q_i^{right}-x_t+x_0\in\{p_j\}$ for all $i = 1, ..., N_o$.
If we denote $q_i^{left}+x_t-x_0:= p^{left}_i$ and $q_i^{right}-x_t+x_0 = p^{right}_i$, the minimization becomes:
\begin{align*}
   \min_{\beta} \bigg\{ \frac{w_1}{N_f} \sum_{i=1}^{N_f} |\mathcal{L}(\mathcal{F}_{\beta}(p_i)) - f(p_i)|^2 + \frac{w_2}{2} \sum_{i=1}^{2} |(\mathcal{F}_{\beta}(q_1) - \mathcal{F}_{\beta}(q_2)|^2 \\
   +  \frac{w_3}{N_o}\sum_{i=1}^{N_0}\bigg( |(\mathcal{F}_{\beta}(q_i^{left}) - \mathcal{F}_{\beta}(p_i^{left})|^2
   +|(\mathcal{F}_{\beta}(q_i^{right}) - \mathcal{F}_{\beta}(p_i^{right})|^2 \bigg) \bigg\},
\end{align*}
where $w_1+w_2+w_3 = 1$ are positive constants.

\section{Numerical examples}
\label{numerical_sec}
In this section, we perform numerical experiments to demonstrate our proposed method (neural homogenization based PINN: NH-PINN). Three equations will be considered; they have the different cell problems and the corresponding homogenized equations are also different.
We will first define all necessary notations in Section (\ref{numerical_notations}). The remaining sections are the experiment results.

It should be noted that our reference solutions are obtained by the finite element method. 
This includes solving the cell problems and the homogenized equations.
To evaluate the homogenized coefficients, one usually needs to calculate the derivatives of the cell problems and do the integration. For the reference solution, we use the finite difference scheme to evaluate the derivatives; 
and we use the quadrature rule to compute the integration.
\subsection{Notations}
\label{numerical_notations}
We first introduce the notations and define the relative errors. Suppose we are solving the following equation with a proper boundary condition:
\begin{align}
    L(a(x, x/\epsilon),u_{\epsilon}) = f, x\in\Omega
    \label{numerical_fine}
\end{align}
where $L$ is a well defined operator; $a(x, x/\epsilon)$ is the coefficient which contributes the multiscale property to the system; and $f$ is the source with a proper regularity. $u_{\epsilon}$ is the reference solution which is obtained by the finite element methods (FEM) with the fine mesh. To get the homogenized equation, we first compute the cell problems; the solutions of the cell problems are the key component in determining the homogenized equation. 
The cell problems can be solved by PINN and the traditional numerical methods.
We denote the the homogenized equation as
\begin{align}
    L^*(b^*(x), v(x)) = f,
    \label{numerical_fem_homog}
\end{align}
if $b^*(x)$ (cell problems) is computed by the FEM; here $v(x)$ is the FEM solution to the system. 
We denote the homogenized equation as
\begin{align}
    L^*(a^*(x), \cdot) = f,
    \label{numerical_pinn_homog}
\end{align}
if $a^*(x)$ (cell problems) is computed by PINN. If the system is solved by the PINN, the solution will be denoted as $p(x)$; and if this equation is solved by the FEM, the solution will be denoted as $w(x)$. 
This solution $w(x)$ can quantify the errors of the PINN method in solving the homogenized equation (\ref{numerical_pinn_homog}). 
The notations of the different solutions are summarized in Table (\ref{table_notations}).
\begin{table}[H]
\centering
\begin{tabular}{||c c c||} 
\hline
Solution notation & Cell problems solver & Homogenized equation solver \\ [0.5ex] 
\hline
 $p(x)$ & PINN & PINN \\ [0.5ex]
\hline
$v(x)$  & FEM & FEM  \\ [0.5ex]
\hline
$w(x)$   & PINN & FEM \\ [0.5ex]
\hline
\end{tabular}
\caption{Notations of the solutions. $u_{\epsilon}$ which is not listed here is the reference solution which solves the multiscale PDE directly by fine scale  finite element methods. }
\label{table_notations}
\end{table}
We compute the relative errors to quantify the accuracy of the proposed method (NH-PINN) and the definitions are given as follow,
\begin{align*}
    &e_1 = \frac{\|p(x) - u_{\epsilon}(x)\|}{\|u_{\epsilon}\|}, e_2 = \frac{\|w(x) - u_{\epsilon}(x)\|}{\|u_{\epsilon}(x)\|}, e_3 = \frac{\|p(x) - w(x)\|}{\|w(x)\|},
    e_4 = \frac{\|v(x) - u_{\epsilon}(x)\|}{\|u_{\epsilon}\|}
\end{align*}
where $\|.\|$ is the $L_2$ norm. 
\textcolor{black}{
$e_1$ is the error of the proposed method (NH-PINN); this error is the ultimately measure of NH-PINN and the number will be compared with the relative error of the classical PINN.  
$e_2$ is the error coming from the cell problems since we fix the homogenized equation solvers; if this error is big, this means the cell problems are not well solved by the PINN.
$e_3$ measures the PINN solver accuracy of solving the homogenized equation; if this error is small,
it can then show that the homogenized equation is easy to be solved by PINN. 
$e_4$ is the error of the homogenization implemented by the classical numerical methods, 
this is theoretical lower bound, NH-PINN cannot be better than this.}

\subsection{2D elliptic equation}
\begin{figure}[H]
\centering
\includegraphics[scale = 0.5]{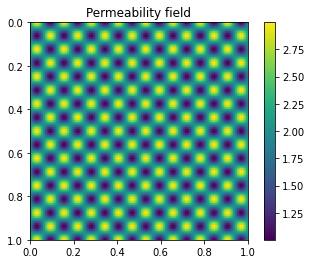}
\includegraphics[scale = 0.5]{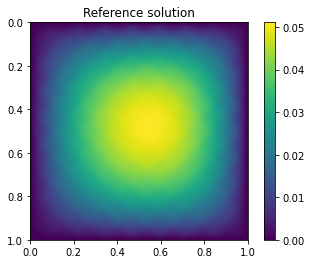}
\caption{Left: permeability of the 2D elliptic problem. Note that $\epsilon = 1/8$. Right: the reference solution.} 
\label{numerical_2d_kappa}
\end{figure}

The first example is a $2D$ elliptic problem; the detailed homogenization is presented in Section (\ref{homog_sec}).
We consider the following $2D$ elliptic equation:
\begin{align}
    -\frac{\partial}{\partial x_i} \big( a(\frac{x}{\epsilon})\frac{\partial }{\partial x_i} u_{\epsilon}(x)\big) = f(x), x\in \Omega,\\
    u_{\epsilon}(x) = 0, x\in\partial \Omega.
    \label{numerical_eqn_2d}
\end{align}
In our examples,  $ \Omega= [0, 1]^2$ and the permeability $a(x/\epsilon) = 2+\sin(2\pi x_1/\epsilon)\cos(2\pi x_2/\epsilon)$ and $\epsilon = \frac{1}{8}$ (shown in Figure (\ref{numerical_2d_kappa}) ). The source is $f(x) = \sin(x_1)+\cos(x_2)$. 

\begin{figure}[H]
\centering
\includegraphics[scale = 0.45]{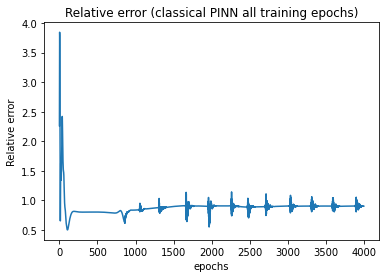}
\includegraphics[scale = 0.45]{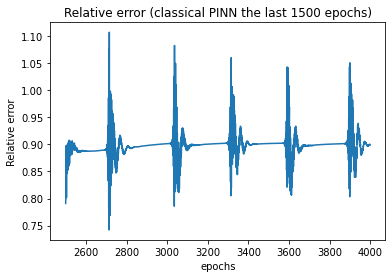}
\caption{$2D$ elliptic problem by the classical PINN. Relative error as a function of the training epochs, the entire history (left), the last $1500$ epochs (right). We train the network for $4000$ epochs with Adam gradient descent (learning rate 0.0001).
The average relative error of the last 500 epochs is $0.900758$; 
} 
\label{pinn_2d_bad_result}
\end{figure}

We first present the results of solving the Equation (\ref{numerical_eqn_2d}) using the classical PINN; we use the same network structure and the result of the prediction is shown in Figure (\ref{pinn_2d_bad_result}).

\begin{figure}[H]
\centering
\includegraphics[scale = 0.45]{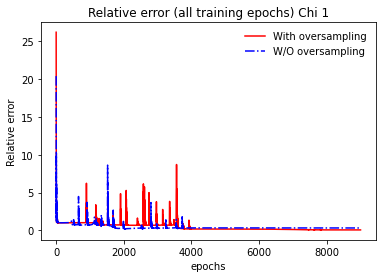}
\includegraphics[scale = 0.45]{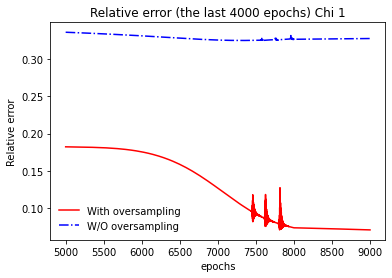}
\caption{$2D$ elliptic cell problem $\chi_1$. Relative error as a function of the training epochs for $\chi_1$. Left: all training epochs; right: the last 4000 epochs.
The average relative errors of the last 300 epochs for the oversampling and without oversampling are $0.072034$ and $ 0.326963$ respectively. } 
\label{numerical_2d_chi1_loss}
\end{figure}

\begin{figure}[H]
\centering
\includegraphics[scale = 0.45]{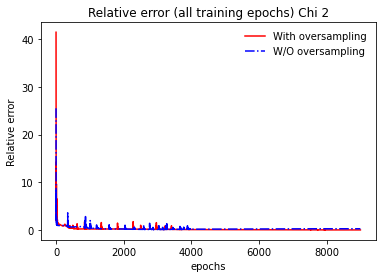}
\includegraphics[scale = 0.45]{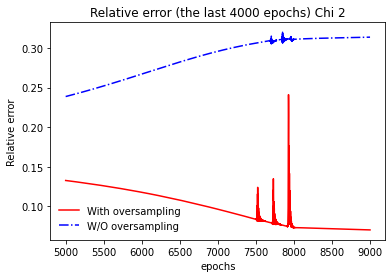}
\caption{$2D$ elliptic cell problem $\chi_2$. Relative error as a function of the training epochs for $\chi_2$. Left: all training epochs; right: the last 4000 epochs.
The average relative errors of the last 300 epochs for the oversampling and without oversampling are $0.071396$ and $0.312961$ respectively.} 
\label{numerical_2d_chi2_loss}
\end{figure}

For NH-PINN, to solve the cell problem, we use a 4-layer network ($2\times 64\rightarrow 64\times 64 \rightarrow 64\times 64 \rightarrow 64\times 1$) activated by Tanh. The domain is discretized with a $101\times 101$ mesh; we use these points and additional 2-layer oversampling in each side of the domain in the training.
We find that the relative errors for the oversampling methods are very robust to the learning rate, the weights of the losses and other hyper-parameters; that is, relative error drops to the similar value even if we change some of the hyper-parameters.
We do not observe the similar results for the standard method without oversampling and hence do a series experiments and then choose the best results. The results for $\chi_1$ and $\chi_2$ are shown in Figure (\ref{numerical_2d_chi1_loss}) and Figure (\ref{numerical_2d_chi2_loss}) respectively.

\begin{figure}[H]
\centering
\includegraphics[scale = 0.45]{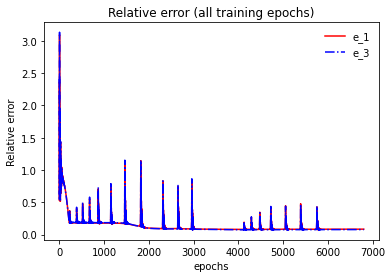}
\includegraphics[scale = 0.45]{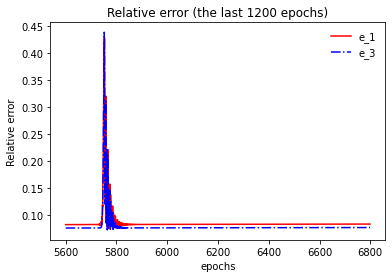}
\caption{$2D$ elliptic cell problem $e_1$ and $e_3$ relative errors with respective to the training epochs. Left: history of all training epochs; right: history of the last 1200 epochs.
The average relative errors of the last 500 epochs are: $e_1 = 0.082994$ and $e_3 = 0.076604$.} 
\label{numerical_2d_e1e3}
\end{figure}

For solving the learnt homogenized equation, we use a network with the same structure as the cell problem. The domain is discretized with a uniform mesh of the size $21\times 21$; and these grid points are then used in the training. To test the network performance, we use a $101\times 101$ mesh.
The relative errors are plotted in Figure (\ref{numerical_2d_e1e3}) and 
we are interested in the relative errors after the training is stabilized, i.e., we compute the average errors of the last $500$ epochs. The numbers are shown in the Table (\ref{numerical_2d_table_relative}). 

\begin{table}[H]
\centering
\begin{tabular}{||c c c c||} 
\hline
$e_1$ & $e_2$ & $e_3$ & $e_4$\\ [0.5ex] 
\hline
$0.082994$ & $0.0212534$  & $0.076604$ & $0.021316$\\ [0.5ex]
\hline
\end{tabular}
\caption{Relative errors for the $2D$ elliptic problem.}
\label{numerical_2d_table_relative}
\end{table}

\begin{figure}[H]
\centering
\includegraphics[scale = 0.45]{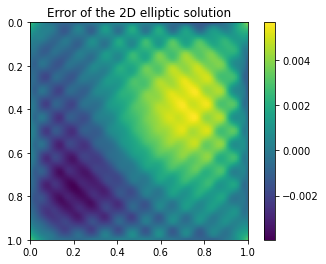}
\caption{$2D$ elliptic error distribution. The error is calculated as the difference between the reference solution and the NH-PINN solution.} 
\label{numerical_2d_diff}
\end{figure}

Since $2D$ solution is not easy to observe, we also calculate the error $u_{reference}- u_{NH-PINN}$ and plot the $2D$ error distribution; please check Figure (\ref{numerical_2d_diff}) for the illustration.

\subsubsection{Interpretation of the results}
We can see from the Table (\ref{numerical_2d_table_relative}) and Figure (\ref{pinn_slow_bad_result}), the $e_1$ drops to around $8\%$; this is a great improvement when compared to the classical method (relative error $0.900758$) which we apply the PINN directly to Equation (\ref{numerical_fine}).
$e_2$ is close to $e_4$ which is the theoretical optimal relative error, this means that the cell problems can be solved by PINN very accurately; besides, $e_1$ is slightly bigger than $e_3$, together these two comparisons imply that
most of the errors come from solving the homogenized equation.
$e_3$ is around $0.076604$, this is the error of solving the homogenized equation using PINN;
however, athe homogenized equation is much easier to be solved by PINN when compared to the multiscale PDE.

One interesting fact is that $e_2<e_4$. We observe the similar phenomenon for some other examples. 
We have an intuitive explanation for this phenomenon. The problem may originate from the evaluation of the $\partial \chi_i/\partial y$.
For the traditional method, the easiest way to calculate the derivative is the finite difference method; computationally, this approximation contributes an error and as a result, 
$a^*$ may not be evaluated exactly.
On the other hand, for NH-PINN method, we approximate $\chi_i$ by the network. The derivatives can be derived by the auto differentiation.
If $\chi_i$ are approximated accurately, we may expect a small error in the derivatives when compared to the exact derivatives. 
Consequently, we obtain a better $a^*$ and this results in a smaller error.
The derivatives can be approximated by some high-order methods; however, the computational costs may be higher than the neural network method.
We observe the similar situation in some other examples, we hence conclude that our proposed NH-PINN method of solving the cell problems potentially can give an accurate approximation of the homogenized coefficients in a relative easier way;
as a result, the final solution of the resulting homogenized equation may also be improved. 

\subsubsection{Transfer learning}
\begin{figure}[H]
\centering
\includegraphics[scale = 0.45]{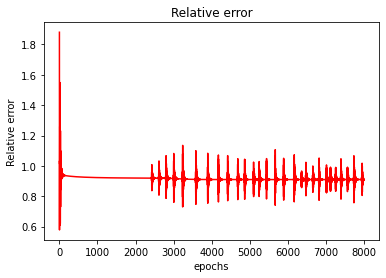}
\caption{2D elliptic transfer learning of classical PINN.  } 
\label{numerical_2d_transfer}
\end{figure}

In this section, we present a transfer learning result which can explain the reason of the PINN failure.
As we have discussed in Section (\ref{pinn_reason_sec}),
we can initialize the training of the classical PINN with a set of parameters which can capture the low frequency component of the solution.
Since we solve the homogenized equation with the network as the classical PINN and the homogenized solution is the coarse scale component, 
the trained network of NH-PINN method can be used the initialization of the classical PINN directly.
All settings including the weights and learning rate are kept the same;
but we get an inaccurate result. 
The relative error blows up very quickly and the stabilised solution has a relative error ($0.909675$, please check Figure (\ref{numerical_2d_transfer}))
which is similar as applying PINN directly to the multiscale PDE. 
\textcolor{black}{This implies that 
the set of parameters which gives an approximation to the coarse scale solution is not optimal (local and global) of the multiscale loss of the classical PINN;
the multiscale loss prevents the network from learning the coarse scale solution. 
}

\subsection{1D slowly varying periodic coefficients}
\label{numerical_sec_1dslow}
\begin{figure}[H]
\centering
\includegraphics[scale = 0.45]{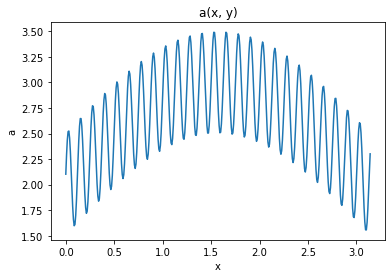}
\includegraphics[scale = 0.45]{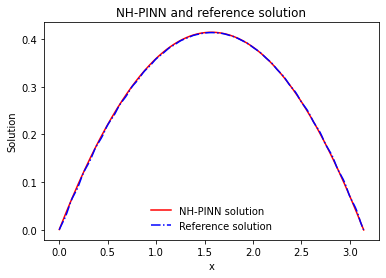}
\caption{1D slowly varying elliptic problem. Left: $a(x, x/\epsilon)$, where $\epsilon = 1/8$. Right: NH-PINN solution vs the reference solution. The relative eror $e_1 = 0.005819$.} 
\label{numerical_1d_slow_axy_sol}
\end{figure}

\begin{figure}[H]
\centering
\includegraphics[scale = 0.3]{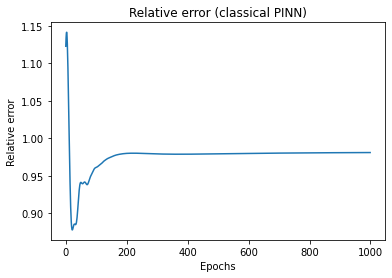}
\includegraphics[scale = 0.3]{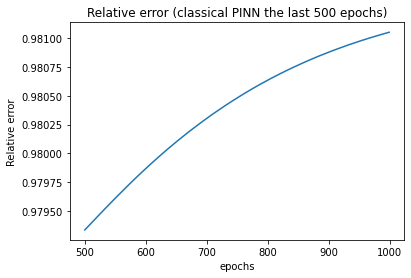}
\includegraphics[scale = 0.3]{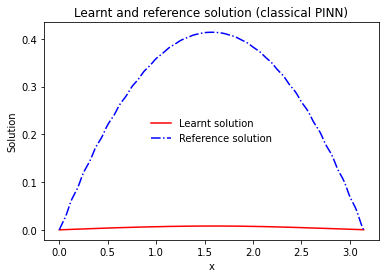}
\caption{$1D$ slowly varying elliptic problem by the classical PINN. Relative error as a function of the training epochs, the entire history(left), the last $500$ epochs (middle) and the solution (right). We train the network for $1000$ epochs with Adam gradient descent (learning rate 0.0001). $\omega_1 = 1/11$ and $\omega_2 = 10/11$. 
The average relative error of the last 100 epochs is $0.980968$; 
however we observe a climb of the error (check the middle image).
We vary the learning rate and the loss weights, however, for all the combinations, we fail to get a satisfactory result.   } 
\label{pinn_slow_bad_result}
\end{figure}

In this section, we consider again an elliptic equation. The difference is that the coefficient $a(x, x/\epsilon)$ has a slow dependency $x$; 
as a result, the homogenization process is different: we need to solve a cell problem for each $x$ in the domain; besides,
the $a^*$ in the homogenized equation also depends on $x$.
We consider the following elliptic equation:
\begin{align}
    -\frac{d}{dx}\big( a(x, \frac{x}{\epsilon}) \frac{d}{dx} u_{\epsilon} \big) = f, x\in [0, \pi],\\
    u_{\epsilon}(0) = u_{\epsilon}(\pi) = 0.
    \label{numerical_eqn_slow}
\end{align}
In this example, $a(x, x/\epsilon) = 0.5\sin(2\pi x/\epsilon)+\sin(x)+2$ where $\epsilon = \frac{1}{8}$, it is demonstrated in Figure (\ref{numerical_1d_slow_axy_sol}) Left, we can observe the high frequency oscillations. The source $f(x) = \sin(x)$.
We first present the results of solving Equation (\ref{numerical_eqn_slow}) using the classical PINN. We use the same network structure as NH-PINN detailed later and the result of the prediction is shown in Figure (\ref{pinn_slow_bad_result}).

\begin{figure}[H]
\centering
\includegraphics[scale = 0.45]{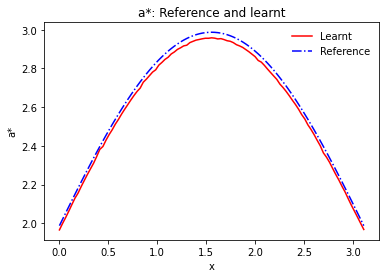}
\includegraphics[scale = 0.45]{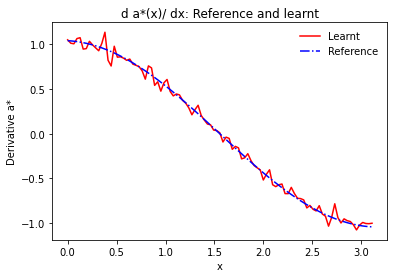}
\caption{$a^*(x)$ and $d a^*(x)/dx$ of the slowly varying problem.} 
\label{numerical_1d_slow_a_star}
\end{figure}

For NH-PINN, to solve each cell problem, we use a network of structure $1\times 64\rightarrow 64\times 64 \rightarrow 64\times 64 \rightarrow 64\times 1$; the network is activated with Tanh as usual. To train the network, we use $101$ uniform spaced points.
We solve the cell problems by PINN and denote the resulting homogenized coefficient as $a^*_p(x)$; the relative error is then calculated as 
$e_a = \|a^*(x)-a^*_p(x)\|/\|a^*(x)\| = 0.010304$. We present $a^*(x)$ and $d a^*(x)/d x$ in Figure (\ref{numerical_1d_slow_a_star}).

\begin{figure}[H]
\centering
\includegraphics[scale = 0.45]{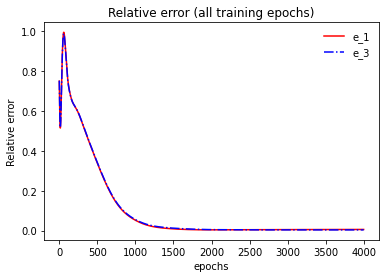}
\caption{1D slowly varying problem $e_1$ and $e_3$ relative errors with respective to the training epochs. 
The average relative errors of the last 500 epochs are: $e_1 = 0.005819$ and $e_3 = 0.003665$.} 
\label{numerical_1d_slow_e1e3}
\end{figure}

Lastly, we need to solve the learnt homogenized equation. We use a network with the same structure as before and $101$ points for training.
The model is tested with $401$ points.
The history of the relative errors $e_1$ and $e_3$ are shown in Figure (\ref{numerical_1d_slow_e1e3}).
We also compute the average relative errors (when the training is stable) in Table (\ref{numerical_1d_slow_table_relative}); more precisely, $e_1$ and $e_3$ are calculated as the average errors of the last $500$ epochs of the training. The solution is demonstrated in Figure (\ref{numerical_1d_slow_axy_sol}) Right. 

\begin{table}[H]
\centering
\begin{tabular}{||c c c c||} 
\hline
$e_1$ & $e_2$ & $e_3$ & $e_4$\\ [0.5ex] 
\hline
$0.005819$ & $0.007668$  & $0.003665$ & $0.006201$\\ [0.5ex]
\hline
\end{tabular}
\caption{Relative errors for the $1D$ slowly varying elliptic problem.}
\label{numerical_1d_slow_table_relative}
\end{table}

\subsubsection{Interpretation of the results}
We can see from the Table (\ref{numerical_1d_slow_table_relative}) and the Figure (\ref{pinn_slow_bad_result}) that,
the relative error $e_1$ of NH-PINN is much better than the one of the classical PINN ($0.980968$). 
This again shows that our method improves the classical PINN performance.
\textcolor{black}{
There is one comment about this example.
Different from the other two examples, the homogenized coefficients $a^*(x)$ can be evaluated exactly. We hence do not observe that $e_2<e_4$ as the other two examples;
however, we can see from the Table (\ref{numerical_1d_slow_table_relative}) that $e_1<e_4$. This is not observed in the other two examples.
This is probably due to the FEM solver accuracy is lower than the PINN.
Since $e_3$ is very small, this indicates that the PINN solver is very accurate. It follows that the total error (error in the homogenized coefficients and error in solving the homogenized equation) of PINN is even smaller than the theoretical upper bound.
}

\subsubsection{Transfer learning}
\begin{figure}[H]
\centering
\includegraphics[scale = 0.45]{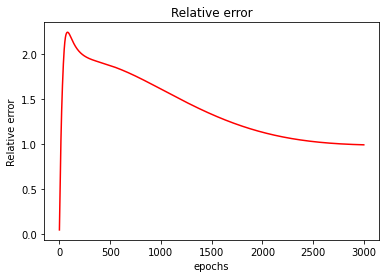}
\caption{1d slow varying elliptic transfer learning of classical PINN.  } 
\label{numerical_1d_transfer}
\end{figure}

Similar as the first example, we use
the trained network of NH-PINN as the initialization of the classical PINN directly.
All settings including the weights and learning rate are kept the same;
but we get an inaccurate result with a relative error $1.005152$.
The relative error increases and stabilizes at the similar value as before (please check Figure (\ref{numerical_1d_transfer})).

\subsection{Diffusion and reaction equation}
\label{numerical_dr_sec}
In this section, we consider another example whose cell problems and the homogenized equation are different with before; The detailed homogenization process can be seen in the Appendix (\ref{app_sec}).
We are going to test the NH-PINN with two different scales. The purpose of the experiments is to show that NH-PINN is able to deal with different scales;
however, the classical PINN fails for both scales and performs even worse when $\epsilon$ is small.
We consider the following example:
\begin{align}
    \frac{\partial u_{\epsilon} }{\partial t}
    -D\nabla\cdot\nabla u_{\epsilon} + \frac{1}{\epsilon}r(\frac{x}{\epsilon})u_{\epsilon} = f, x\in \Omega,\\
    u_{\epsilon}(x) = 0, x\in\partial \Omega.
    \label{numerical_eqn_rd}
\end{align}
In our example, $\Omega = [-\pi, \pi]$. We set $r(x/\epsilon) = \cos(x/\epsilon)$. We use different $\epsilon$ and present the results in the following two sections.
\subsubsection{Scale $\epsilon = 1/10$}
\begin{figure}[H]
\centering
\includegraphics[scale = 0.45]{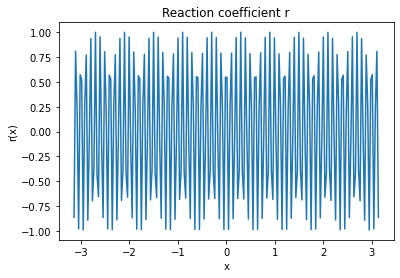}
\includegraphics[scale = 0.45]{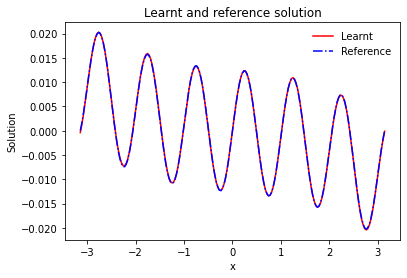}
\caption{Diffusion and reaction problem. Left: $r(x) = cos(x/\epsilon)$, where $\epsilon = 1/10$. Right: NH-PINN solution vs the reference solution. The relative error $e_1 = 0.012388$.} 
\label{numerical_dr_r_sol}
\end{figure}

\begin{figure}[H]
\centering
\includegraphics[scale = 0.3]{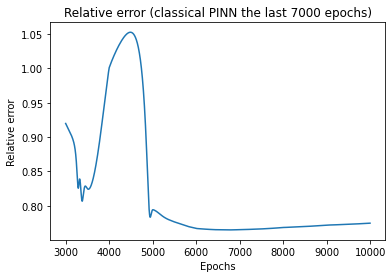}
\includegraphics[scale = 0.3]{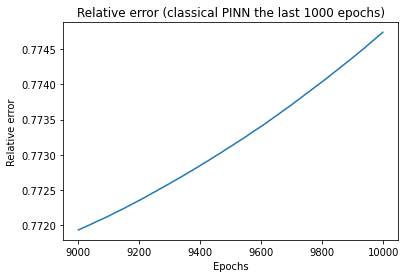}
\includegraphics[scale = 0.3]{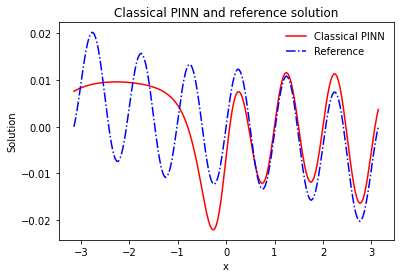}
\caption{Diffusion and reaction problem with $\epsilon = 1/10$ solved by the classical PINN. Relative error of the last $7000$ epochs (left), the last $1000$ epochs (middle) and the solution (right) of the classical PINN. We train the network for $10000$ epochs with Adam gradient descent (learning rate 0.015). $\omega_1 = 1/7$, $\omega_2 = 5/7$ and $\omega_3 = 1/7$, where $\omega_3$ is the weight of the initial condition. The average relative error of the last 100 epochs is $0.774558$; however we observe a climb of the error (check the middle image). We vary the learning rate and the loss weights, however, for all the combinations, we fail to get a satisfactory result.   } 
\label{pinn_dr_result}
\end{figure}

In this section, we set $\epsilon = \frac{1}{10}$ and $D = 2$, the $r(x/\epsilon)$ is demonstrated in Figure (\ref{numerical_dr_r_sol}) Right. The source is $f(x) = \sin(2\pi x)$. 
As we have discussed before, the classical PINN cannot give us an accurate prediction; the results of applying PINN directly are shown in Figure (\ref{pinn_dr_result}).

\begin{figure}[H]
\centering
\includegraphics[scale = 0.45]{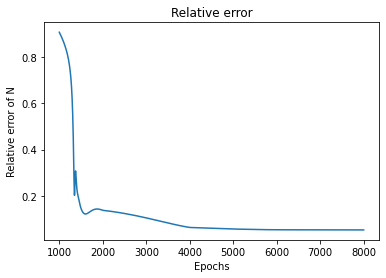}
\includegraphics[scale = 0.45]{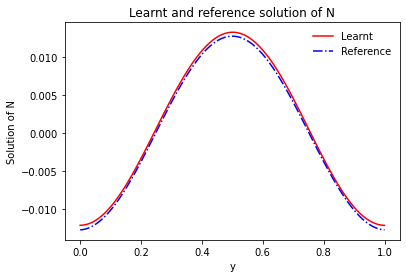}
\caption{Diffusion and reaction equation cell problem (\ref{appendix_n}) with $\epsilon = 1/10$. Left: relative error with respected to the training epochs; note: only last 7000 epochs are shown. The average relative error in $N(y)$ of the last $500$ epochs is $0.053309$. Right: N(y) reference vs learnt.} 
\label{numerical_rd_N}
\end{figure}

For NH-PINN, the cell problem is solved with a 4-layer network of the structure $1\times 64\rightarrow 64\times 64 \rightarrow 64\times 64 \rightarrow 64\times 1$;
the network is activated by the Tanh function. For the network training, we use $101$ points uniformed placed in the domain and additional 2-layer oversampling on each side of the domain. The relative error in $N$ drops to $0.053309$ when the training is stable (we take the average of the last 500 epochs).
The solution and training history of the cell problem are shown in Figure (\ref{numerical_rd_N}).

\begin{figure}[H]
\centering
\includegraphics[scale = 0.45]{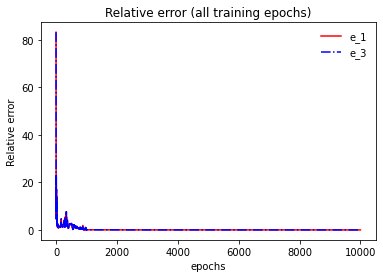}
\includegraphics[scale = 0.45]{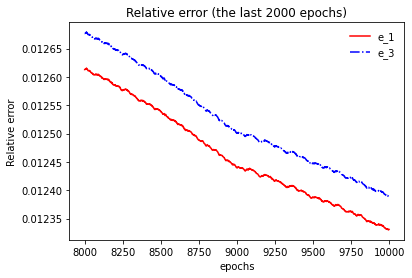}
\caption{Diffusion and reaction problem with $\epsilon =1/10$, $e_1$ and $e_3$ relative error as a function of the training epochs. Left: history of all training epochs; right: history of the last 2000 epochs.
The average relative errors of the last 500 epochs are: $e_1 = 0.012388$ and $e_3 = 0.012448$.} 
\label{numerical_rd_e1e3}
\end{figure}

To solve the homogenized equation, we use another 4-layer network activated by Tanh; the network structure is the same to the one used in the cell problems except that the first layer's dimension is $2\times 64$. The spatial domain and the temporal domain are both discretized with $100$ uniformed placed mesh points; the grid points are then used in the training.
We will test the the solution at the terminal time and uniformly place $201$ points in space.
The training epoch is set to be $10,000$ and we use the standard Adam gradient descent algorithm.
The history of the relative errors are shown in Figure (\ref{numerical_rd_e1e3}); and the average relative errors when the training is stable are presented in Table (\ref{numerical_1d_dcr_table_relative}). We also solve Equation (\ref{numerical_eqn_rd}) using the classical PINN; we use the same network structure and the result of the prediction is shown in Figure (\ref{pinn_dr_result}).

\begin{table}[H]
\centering
\begin{tabular}{||c c c c||} 
\hline
$e_1$ & $e_2$ & $e_3$ & $e_4$\\ [0.5ex] 
\hline
$0.012388$ & $0.0012916$  & $0.012448$ & $0.0012917$\\ [0.5ex]
\hline
\end{tabular}
\caption{Relative errors for the diffusion reaction problem.}
\label{numerical_1d_dcr_table_relative}
\end{table}

\begin{figure}[H]
\centering
\includegraphics[scale = 0.45]{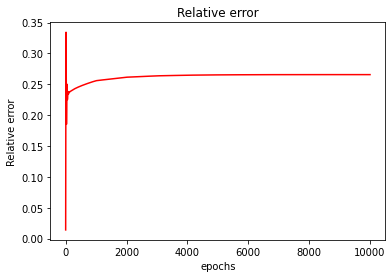}
\caption{Diffusion reaction equation with $\epsilon = 1/10$ transfer learning of classical PINN.  } 
\label{numerical_diffusion_transfer}
\end{figure}

From the Table (\ref{numerical_1d_dcr_table_relative}) and the Figure (\ref{pinn_dr_result}),
firstly we can observe that the relative error $e_1$ of NH-PINN is much better than the classical PINN.
If we look at the Figure (\ref{numerical_rd_e1e3}), the relative errors of NH-PINN are still decaying.
Since $e_2$ and $e_4$ are very closed to each other ($e_2<e_4$),
we conclude that the PINN aided homogenization is a potential alternative to the traditional numerical driven homogenization.
Finally, because $e_2$ is small, this implies that most errors of the method still come from solving the homogenized equation with PINN; 
however, this has been much improved when compared to applying PINN on the multiscale PDE. The transfer learning is the same as before, the stabilized relative error is $0.265721$  
and cannot give us a better result (please check Figure (\ref{numerical_diffusion_transfer})).

\subsubsection{Scale $\epsilon = 1/50$}
\begin{figure}[H]
\centering
\includegraphics[scale = 0.45]{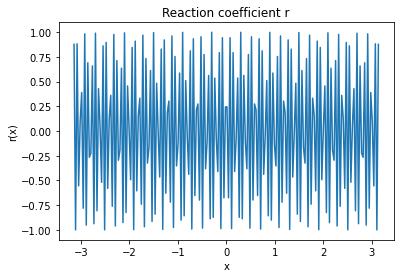}
\includegraphics[scale = 0.45]{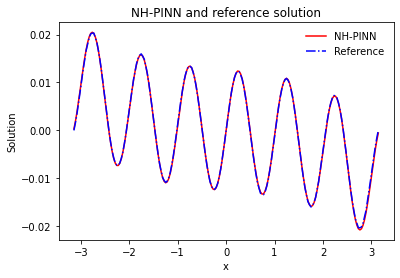}
\caption{Diffusion and reaction problem. Left: $r(x) = cos(x/\epsilon)$, where $\epsilon = 1/50$. Right: NH-PINN solution vs the reference solution. The relative error $e_1 = 0.018582$.} 
\label{numerical_dr_r_sol_50}
\end{figure}

\begin{figure}[H]
\centering
\includegraphics[scale = 0.3]{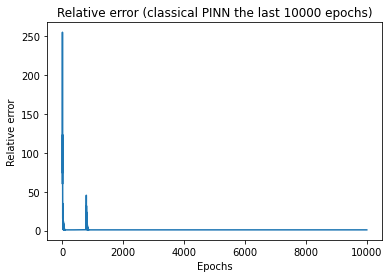}
\includegraphics[scale = 0.3]{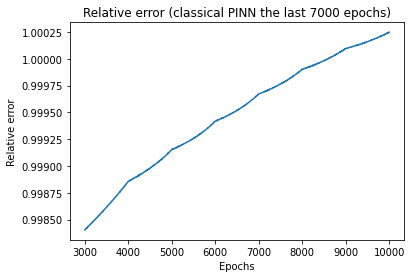}
\includegraphics[scale = 0.3]{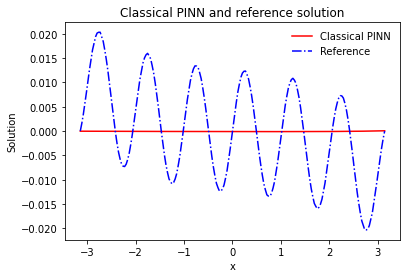}
\caption{Diffusion and reaction problem with $\epsilon = 1/50$ solved by the classical PINN. Relative error of the last $10,000$ epochs (left), the last $7,000$ epochs (middle) and the solution (right) of the classical PINN. We train the network for $10000$ epochs with Adam gradient descent (learning rate 0.05). $\omega_1 = 1/7$, $\omega_2 = 5/7$ and $\omega_3 = 1/7$, where $\omega_3$ is the weight of the initial condition. The average relative error of the last 100 epochs is $1.000237$; however we observe a climb of the error (check the middle image). We vary the learning rate and the loss weights, however, for all the combinations, we fail to get a satisfactory result.   } 
\label{pinn_dr_result_50}
\end{figure}

In this section, we set $\epsilon = \frac{1}{50}$ and $D = 2$, the $r(x/\epsilon)$ is demonstrated in Figure (\ref{numerical_dr_r_sol_50}) Right. The source is $f(x) = \sin(2\pi x)$. 
As we have discussed before, the classical PINN cannot give us an accurate prediction; the results of applying PINN directly are shown in Figure (\ref{pinn_dr_result_50}).

\begin{figure}[H]
\centering
\includegraphics[scale = 0.45]{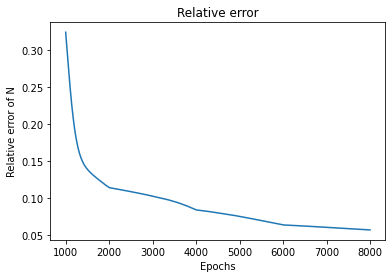}
\includegraphics[scale = 0.45]{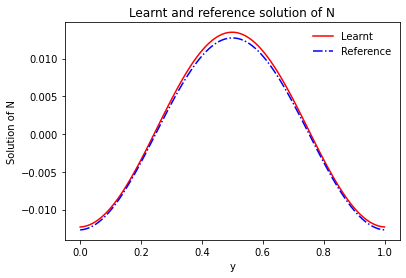}
\caption{Diffusion and reaction equation cell problem (\ref{appendix_n}) with $\epsilon = 1/50$. Left: relative error with respected to the training epochs; note: only last 7000 epochs are shown. The average relative error in $N(y)$ of the last $500$ epochs is $0.057602$. Right: N(y) reference vs learnt.} 
\label{numerical_rd_N_50}
\end{figure}

For $\epsilon = 1/50$, we use the same network setting as before. The relative error of $N$ drops to $0.057602$ when the training is stable (we take the average of the last 500 epochs).
The solution and training history of the cell problem are shown in Figure (\ref{numerical_rd_N_50}).

\begin{figure}[H]
\centering
\includegraphics[scale = 0.45]{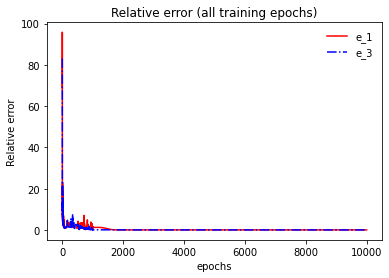}
\includegraphics[scale = 0.45]{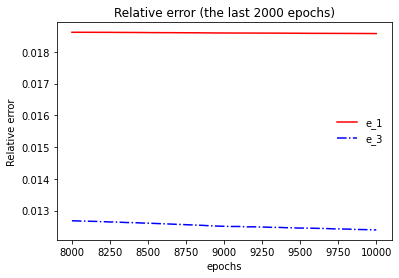}
\caption{Diffusion and reaction problem 
with $\epsilon = 1/50$, $e_1$ and $e_3$ relative error as a function of the training epochs. Left: history of all training epochs; right: history of the last 2000 epochs.
The average relative errors of the last 500 epochs are: $e_1 = 0.012388$ and $e_3 = 0.022184$.} 
\label{numerical_rd_e1e3_50}
\end{figure}

For the homogenized equation with $\epsilon = 1/50$, we use the same setting as the $\epsilon = 1/10$ case.
The history of the relative errors are shown in Figure (\ref{numerical_rd_e1e3_50}); and the average relative errors when the training is stable are presented in Table (\ref{numerical_1d_dcr_table_relative_50}). We also solve Equation (\ref{numerical_eqn_rd}) using the classical PINN; we use the same network structure and the result of the prediction is shown in Figure (\ref{pinn_dr_result_50}).

\begin{table}[H]
\centering
\begin{tabular}{||c c c c||} 
\hline
$e_1$ & $e_2$ & $e_3$ & $e_4$\\ [0.5ex] 
\hline
$0.018582$ & $0.011637$  & $0.022184$ & $0.0116435$\\ [0.5ex]
\hline
\end{tabular}
\caption{Relative errors for the diffusion reaction problem, $\epsilon = 1/50$.}
\label{numerical_1d_dcr_table_relative_50}
\end{table}

\begin{figure}[H]
\centering
\includegraphics[scale = 0.45]{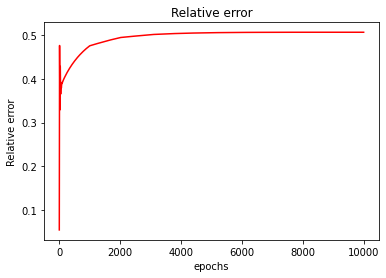}
\caption{Diffusion reaction equation with $\epsilon = 1/50$ transfer learning of classical PINN.  } 
\label{numerical_diffusion_transfer_50}
\end{figure}

Similar as the $\epsilon = 1/10$ case, 
firstly we can observe from the Table (\ref{numerical_1d_dcr_table_relative_50}) and the Figure (\ref{pinn_dr_result_50}) that the relative error $e_1$ of NH-PINN is much better than the classical PINN.
Since $e_2<e_4$ are closed to each other,
we conclude that the PINN aided homogenization is a potential alternative to the traditional numerical driven homogenization.
The transfer learning (check Figure (\ref{numerical_diffusion_transfer_50}) ) is the same 
as before ( with a relative error $0.506884$)
and cannot give us a better result.

\section{Conclusion}
Multiscale problems widely exist in the real life; to solve multiscale problems, fine scale solvers are required, but the computational cost is very high.
Researchers consider solving the multiscale problem with the data-driven approaches.
In this work, we propose to solve the multiscale problems by the physics-informed neural network (PINN) with the help of the homogenization.
We first find that the classical PINN is unable to solve the multiscale problems. Homogenization is a 3-step PDE technique which is used to approximate the multiscale equation by a homogenized equation. We find that all steps of the homogenization can be implemented by the PINN accurately.
In particular, we propose an oversampling strategy which greatly improves PINN accuracy for solving the periodic problems; this technique is used in the first step of the homogenization. We conduct several experiments and find that the accuracy of PINN is greatly improved by our method. We also observe that PINN assisted homogenization is also an accurate approach of implementing homogenization; we hence conclude that the PINN is a potential alternative to the numerical driven homogenization.  


\appendix

\section{Homogenization of the diffusion reaction equation}
\label{app_sec}
In this section, we provide the details of the homogenizing of a diffusion reaction equation. 
The problem is defined as:
\begin{align}
    \frac{\partial u_{\epsilon}}{\partial t} -D\nabla\cdot \nabla u_{\epsilon}+\frac{1}{\epsilon}r(x/\epsilon) u_{\epsilon} = f,
\end{align}
where $\int_Y r(y) dy = 0$ for the solvability; $Y = [0, 1]^d$ is the unit cube.
We seek the same asymptotic expansion as (\ref{homg_asym}) and by equating the power of $\epsilon$, we have:
\begin{align*}
    -D\nabla\cdot \nabla u_0 = 0
\end{align*}
We have $u_0$ is independent of $y$ and we can then further simplifies the $\epsilon^{-1}$ term, 
\begin{align*}
    -D\nabla\cdot \nabla u_1 = r_0u_0,
\end{align*}
where $u_1(x, y)$ is double periodic in $y$ with period $Y$.
The cell problem can then be defined as:
\begin{align}
    -D\nabla\cdot \nabla N(y) = -r(y), y\in Y.
    \label{appendix_n}
\end{align}
The problem has the double periodic boundary condition and we assume $N(y)$ has zero average in $Y$;
the $\int_Y r(y) dy = 0$ guarantees the equation is solvable.
$u_1$ can be then expressed as:
\begin{align*}
    u_1(x, y) = N(y)u_0(x).
\end{align*}
Finally for $\epsilon^0$ term,
\begin{align*}
    \frac{\partial u_0}{\partial t}-D\nabla\cdot \nabla u_0- D\Delta_{xy}u_1 - D\nabla\cdot \nabla u_2+r(t)u_1 = f(x).
\end{align*}
Integrate the above equation over $y$ in one period $Y$, we finally have,
\begin{align}
     \frac{\partial u_0}{\partial t} - D\nabla\cdot \nabla u_0+r^* u_0 = f,
\end{align}
where the homogenized coefficient $r^*$ is:
\begin{align*}
    \int_Y r(y)N(y) dy.
\end{align*}

\bibliographystyle{abbrv}
\bibliography{references}
\end{document}